\documentclass[12pt]{amsart}
\usepackage{amsmath}
\usepackage{amssymb}
\usepackage{verbatim}
\usepackage{graphicx}
\usepackage{color}
\usepackage{dsfont}
\usepackage{hyperref}
\hypersetup{colorlinks=true, linkcolor=blue, citecolor=magenta, urlcolor=wine}
\usepackage{url}
\usepackage{fancyhdr}
\usepackage{tikz-cd}
\usepackage{tikz}
\usepackage{lipsum}
\usepackage{adjustbox}
\usetikzlibrary{shapes.geometric,arrows.meta}
\tikzcdset{arrow style=tikz}
\usepackage{calc}
\usepackage{stmaryrd}
\usepackage{tikz-cd}
\usepackage{lipsum}
\usepackage{adjustbox}
\usepackage{quiver}

\newcommand{\arxiv}[1]{\href{http://arxiv.org/abs/#1}{\tt arXiv:\nolinkurl{#1}}}
\newcommand{\arXiv}[1]{\href{http://arxiv.org/abs/#1}{\tt arXiv:\nolinkurl{#1}}}
\newcommand{\doi}[1]{\href{http://dx.doi.org/#1}{{\tt DOI:#1}}}

\newcommand{\googlebooks}[1]{(preview at \href{http://books.google.com/books?id=#1}{google books})}

\usepackage[new]{old-arrows}  
\tikzstyle{none}=[]
\pgfdeclarelayer{nodelayer}
\pgfdeclarelayer{edgelayer}
\pgfsetlayers{edgelayer,main,nodelayer}

\usepackage{array}
\usepackage{enumerate}
\usepackage{xfrac}
\usepackage{mathrsfs}   
\usepackage{wasysym}	
\usepackage{xcolor}

\usepackage{mathtools}

\usepackage{hyperref}
\hypersetup{
    colorlinks=true,
    linkcolor=blue,
    filecolor=purple,      
    urlcolor=[rgb]{0 0 .6},
}

\newcommand{\1}{\mathds{1}}

\newcommand{\Hom}{\mathrm{Hom}}

\newcommand{\Aut}{\mathrm{Aut}}

\renewcommand{\Vec}{\mathrm{Vec}}

\newcommand{\Rep}{\mathrm{Rep}}

\newcommand{\Id}{\mathrm{Id}}
\newcommand{\Triv}{\mathrm{Triv}}

\renewcommand{\lim}{\mathop{\mathrm{lim}}}

\makeatletter
\g@addto@macro\th@plain{\thm@headpunct{}}
\g@addto@macro\th@definition{\thm@headpunct{}}
\g@addto@macro\th@remark{\thm@headpunct{}}
\makeatother

\makeatletter
\newtheorem*{rep@theorem}{\rep@title}
\newcommand{\newreptheorem}[2]{%
\newenvironment{rep#1}[1]{%
 \def\rep@title{#2 \ref{##1}}%
 \begin{rep@theorem}}%
 {\end{rep@theorem}}}
\makeatother

\newtheorem{theorem}{Theorem}[section]
\newreptheorem{theorem}{Theorem}
\newtheorem{proposition}[theorem]{Proposition}

\newtheorem{lemma}[theorem]{Lemma}

\newtheorem{convention}[theorem]{Convention}

\theoremstyle{definition}
\newtheorem{definition}[theorem]{Definition}

\tikzset{Rightarrow/.style={double equal sign distance,>={Implies},->},
triple/.style={-,preaction={draw,Rightarrow}},
quadruple/.style={preaction={draw,Rightarrow,shorten >=0pt},shorten >=1pt,-,double,double
distance=0.2pt}}

\title{Equivariant fixed points on Tambara-Yamagami categories}
\author{Alexander Betz, David Green and Magnus Hanson}
\date{}

\begin{document}

\begin{abstract}
We specialize recent work of Galindo and Jones classifying equivariant fusion subcategories to the Tambara-Yamagami categories. In this case, we show how their parameterization of subcategories can be reduced to group theoretic data, then provide some explicit computations of the number and type of these subcategories when the acting group is $S_4$. 
\end{abstract}

\maketitle
\tableofcontents

\section{Introduction}

Fusion categories are intricate mathematical objects that simultaneously generalize finite groups and their categories of representations. Classical examples of fusion categories include $\Vec$, the category of finite-dimensional vector spaces, $\Vec_A$ the category of $A$-graded vector spaces and $\Rep(G)$ the category of finite-dimensional representations of a finite group $G$ (over a field of characteristic zero). However, fusion categories go far beyond categories of structured vector spaces - see for instance the famous Fibonacci category \cite{MR1981895}. 

 These categories appear as invariants of subfactors \cite{MR1645078}, descriptions of topologically ordered systems, and in quantum computing \cite{MR2640343}. The utility of fusion categories has encouraged their study as fundamental objects \cite{MR2183279}. Within this study, one research program concerns the construction of new fusion categories from a given one. 

This leads us to study the equivariantization of a group action on a fusion category, a categorical analogue of a group of fixed points. Indeed, given a group $G$ acting by monoid homomorphisms on a monoid $M$, there is a natural notion of the monoid of fixed points $M^G$. This can be regarded as a construction of a new monoid $M^G$ from the data of a given monoid $M$, together with additional data (the action of $G$ on $M$). Categorifying, finite groups act on a monoidal category via a monoidal functor $\underline{G} \to \mathrm{Aut}_\otimes(\mathcal{C})$. For a fusion category $\mathcal{C}$ equipped with such an action, there is a categorified notion of the fixed points, called the $G$-equivariantization and denoted $\mathcal{C}^G$. When $\mathcal{C}$ is fusion and $G$ is finite,  $\mathcal{C}^G$ is also a fusion category and so may have fusion subcategories. These were classified in terms of \textit{equivariant trivializations} by Galindo and Jones \cite{Galindo2025}. 

We will apply their machinery to the famous Tambara-Yamagami categories \cite{MR1659954} which arise from a finite abelian group $A$ equipped with a symmetric nondegenerate bicharacter $\chi$ and a sign choice $\tau = \pm 1/\sqrt{|A|}$. These categories appear regularly as examples in the theory of fusion categories and topological phases \cite{thorngrenwangTY}, and are well understood via the lens of extension theory. 
Through this lens, it is possible to parameterize $G$-actions on Tambara-Yamagami categories in terms of the $\chi$-preserving automorphisms of $A$ and cohomological data \cite{ECM22, galindo2011clifford}.

\begin{reptheorem}{thm:MainTheorem}
    Let $\rho \colon G \to O(A, \chi)$ be a homomorphism of groups, and $[\varphi] \in H^2(G, \mathbb{Z}/2\mathbb{Z})$ provide a $G$-action on $\mathcal{TY}(A, \chi, \tau)$. Then fusion subcategories of $\mathcal{TY}(A, \chi, \tau)^G$ occur in two families, which are respectively parameterized by 
    \begin{itemize}
        \item Triples $(H, B, \eta)$ such that 
        \begin{itemize}
            \item $H \lhd G$
            \item $B \leq A$ where $B$ is $G$-invariant and pointwise fixed by $H$
            \item $\eta \in \Hom(H, \Hom(B, \mathbb{K}^\times))$ a fixed point of the natural $G$ action -  i.e, satisfying $$\eta^x_{g^{-1}hg}=\eta_h^{g_*(x)}$$ 
        \end{itemize}
        \item A normal subgroup $H \lhd G$ such that $H\leq \ker\rho$ with a trivialization $\eta$ of $[\varphi|_H]$ satisfying $$\eta(g^{-1}hg)\varphi(h,g) = \eta(h)\varphi(g,g^{-1}hg)$$
    \end{itemize}
\end{reptheorem}

The article is primarily devoted to the proof of Theorem \ref{thm:MainTheorem}. Section \ref{sec:Examples} is devoted to making Theorem \ref{thm:MainTheorem} explicit for $G=S_4$ and $A=\mathbb{Z}_6 \times \mathbb{Z}_6$, for particular choices of the homomorphism $\rho$ and bicharacter $\chi$ and arbitrary choice of cohomology class $[\varphi]$. The final result is presented in a series of tables.

\textbf{AI Use Statement:} The results presented in Theorem \ref{thm:Example} were initially computationally derived by MH via the use of AI-written code. In the process of removing some of this dependence on machine computation, ChatGPT 6.2 Astra was prompted by DG to provide an effective procedure in terms of exact sequences to compute fixed points.  The resulting output was digested, verified, attributed, and rewritten by all authors and became Lemma \ref{lemma:EffectiveCohomology}. ChatGPT also provided feedback on a draft of this manuscript. All mathematical
statements, proofs, computational results, and references were independently checked by the authors, who take sole responsibility for the contents. 

\textbf{Contributions and Acknowledgments:} This article is the result of MH's undergraduate research project, which was supervised primarily by the other two authors. The majority of the mathematical content, including derivation of the main results and the computations of the examples, primarily originated from MH. All the authors would like to thank Corey Jones for suggesting this classification problem, as well as for his continued mentorship. Access to ChatGPT for DG and MH was supported through the ChatGPT for Academic Researchers program by OpenAI.

\section{Preliminaries}

 A \textit{fusion category} is a finite semisimple $\mathbb{K}$-linear monoidal category with a simple unit object $1$. Here $\mathbb{K}$ is assumed to be an algebraically closed field of characteristic 0. For convenience, we assume all fusion categories are skeletal, with only one object per isomorphism class. We refer the reader to \cite{MR3242743, MR4806154} for additional background on fusion categories, and continue to recall the basic notions of group actions on fusion categories. 
 
The monoidal category $\underline{G}$ has objects indexed by the set $G$ and only identity morphisms. The tensor product is  $g\otimes h\cong gh$ and the associator is the identity. An \textit{action} of $G$ on a fusion category is a monoidal functor $(F,J)$, $F:\underline{G}\rightarrow \Aut_\otimes(C)$, where $\Aut_\otimes(C)$ is the monoidal category of monoidal autoequivalences of $\mathcal{C}$. This assigns a group element $g$ to a monoidal autoequivalence $g_*$ and provides a natural isomorphism $J_{g,h}:g_*\circ h_*\rightarrow (gh)_*$ satisfying 
\begin{center}

\begin{tikzcd}[column sep=huge, row sep=large]
g_* h_* k_*
  \arrow[r, "J_{g,h}\,\mathrm{id}_{k_*}"]
  \arrow[d, "\mathrm{id}_{g_*}\,J_{h,k}"']
& (gh)_* k_*
  \arrow[d, "J_{gh,k}"] \\
g_* (hk)_*
  \arrow[r, "J_{g,hk}"']
& {(ghk)}_*
\end{tikzcd}
\end{center}

\begin{definition}\label{defn:Trivialization}
    A \textit{trivialization} of a $G$-action $(F,J)$ is a monoidal natural isomorphism $\eta:(F,J)\Rightarrow (\Triv,\Triv)$ to the trivial monoidal functor. This corresponds to the following data:

\begin{itemize}
    \item For each object $g$ in $\underline{G}$, a natural isomorphism $\eta_g:g_*\rightarrow\Id$, satisfying the monoidal natural transformation axiom \[\eta_{g}^xg_*(\eta_h^x)=\eta_{gh}^xJ_{g,h}^x.\]
\end{itemize}
\end{definition}

\begin{definition}\cite{MR3242743}
    A $G$-equivariant object is an object $(X,u)$ where $X$ is an object in $\mathcal{C}$ and $u:g_*(X)\rightarrow X$ is a natural isomorphism satisfying,
\begin{center}
\begin{tikzcd}[column sep=large, row sep=large]
g_*\bigl(h_*(X)\bigr)
  \arrow[r, "g_*(u_h)"]
  \arrow[d, "J_{g,h,X}"']
& g_*(X)
  \arrow[d, "u_g"] \\
gh_*(X)
  \arrow[r, "u_{gh}"']
& X
\end{tikzcd}
\end{center}
\end{definition}

together with a unit condition. The category of $G$-equivariant objects, together with 1-morphisms commuting with the $u_g$, assemble into a category $\mathcal{C}^G$ of $G$-\textit{equivariant objects}.  When $\mathcal{C}$ is a fusion category $\mathcal{C}^G$ is also a fusion category \cite{MR3059899}. The $G$-equivariant structure on a tensor product is built from those of the tensorands and the coherence morphisms of the action. 

The results of \cite{Galindo2025} outline how to use trivializations to find all fusion subcategories of the equivariantization $\mathcal{C}^G$.

\begin{theorem}\label{thm:GJ subcats}\cite[Theorem  4.8]{Galindo2025} 
    Let $\mathcal{C}$ be a fusion category equipped with a $G$-action. Fusion subcategories of $\mathcal{C}^G$ are parametrized by triples $(\mathcal{E},H,\eta)$ consisting of:

    \begin{enumerate}
        \item $\mathcal{E}\subset \mathcal{C}$ a $G$-invariant fusion subcategory
        \item $H \trianglelefteq G$ is a normal subgroup
        \item $\eta$ is a $G$-equivariant trivialization of $H$ on $\mathcal{E}$.
    \end{enumerate}
\end{theorem}

We specialize this theorem in Section \ref{sec:Theory} when the category $\mathcal{C}$ is the Tambara-Yamagami fusion category $\mathcal{TY}(A,\chi,\tau)$. There we provide a description of the set of equivariant trivializations purely in terms of the underlying group theoretical and cohomological data of the category $\mathcal{TY}(A,\chi,\tau)$ and $G$-action. To provide such a description, we need to recall the structure of and notation for the categories $\mathcal{TY}(A,\chi,\tau)$. 

Let $A$ be a finite abelian group. The fusion ring $\mathsf{TY}(A)$ has $\mathbb{Z}_+$ basis consisting of the group elements $a\in A$ and an additional element $m$. The multiplication is defined on $m$ and for all $a, b \in A$ by the rules:
$$a \otimes b = ab, \quad a \otimes m = m = m \otimes a, \quad m \otimes m = \bigoplus_{a \in A} a$$
We notice that this fusion ring is $\mathbb{Z}_2$ graded where $A$ has the trivial grading and $m$ has the $1$ grading since $m\otimes m$ is contained in the trivial graded component $A$. This implies that the fusion category $\mathcal{TY}(A,\chi,\tau)$ is $\mathbb{Z}_2$ graded and the universal grading group is $\mathbb{Z}_2$.

\begin{theorem}[{\cite[Theorem 3.2]{MR1659954}}]
Let $A$ be a finite abelian group, let
$\tau=\pm1/\sqrt{|A|}$,
and let $\chi\colon A\times A\to \mathbb{K}^{\times }$ be a symmetric nondegenerate
bicharacter.
We define a fusion category $\mathcal{C}(A,\chi,\tau)$ by taking
the underlying fusion ring to be $\mathsf{TY}(A)$, the unitor isomorphisms
to be identity, and the associators to be
\begin{align*}
    \alpha_{a,b,c} &= 1_{abc}, \\
    \alpha_{a,b,m} = \alpha_{m,a,b} &= 1_{m}, \\
    \alpha_{a,m,b} &= \chi(a,b)\cdot 1_{m}, \\
    \alpha_{a,m,m} = \alpha_{m,m,a} &= \bigoplus_{b\in A} 1_{b}, \\
    \alpha_{m,a,m} &= \bigoplus_{b\in A} \chi(a,b)\cdot 1_b, \\
    \alpha_{m,m,m} &= (\tau\chi(a,b)^{-1}\cdot 1_m)_{a,b}.
\end{align*}
\end{theorem}

We will use $G$ interchangeably for both a group $G$ and the monoidal 1-category whose objects are $G$ and whose morphisms are all identities. This category is \textit{not} the same as $\mathbf{B}G$. By definition, group actions on the Tambara-Yamagami categories are the same thing as monoidal functors $G \to \Aut_\otimes(\mathcal{C}(A,\chi,\tau))$. If we want to understand these functors, we need to understand the target, which can be described in terms of $O(A, \chi)$ - the automorphisms of $A$ preserving $\chi$. 

\begin{proposition} As monoidal categories,
$$\Aut_{\otimes}(\mathcal{TY}(A, \chi, \tau)) \cong O(A, \chi) \times B\mathbb{Z}_2$$
\end{proposition}
\begin{proof}
   The group of objects of $\Aut_{\otimes}(\mathcal{TY}(A, \chi, \tau))$ is identified by several authors as $O(A, \chi)$. See for instance \cite{tambara2000representations, nikshych2008non, ECM22}. The automorphisms of the identity are classified by \cite{GALINDO_2022} as $B\mathbb{Z}_2$. A direct calculation using the representatives provided by \cite{ECM22} shows that this 2-group is split. 
\end{proof}

Now we can start to understand all functors: 

\begin{proposition}\label{prop:monfun G into Aut classification}
    The set of equivalence classes of monoidal functors from $G$ to $\text{Aut}_{\otimes}(\mathcal{TY}(A, \chi, \tau))$ is parameterized by pairs $(\rho, [\varphi])$, where:

\begin{itemize}
    \item $\rho: G \to O(A, \chi)$ is a group homomorphism.
    \item $[\varphi] \in H^2(G, \mathbb{Z}_2)$ is a second cohomology class with coefficients in the trivial $G$-module $\mathbb{Z}_2$.
\end{itemize} 
\end{proposition}
\begin{proof}
    This follows from Theorem 5.5 of \cite{galindo2011clifford}. 
\end{proof}

\begin{convention} Let $G$ act on $\mathcal{TY}(A,\chi,\tau)$. 
     By \cite{ECM22}, we will assume that $\mu_g$ is an identity for all $g$; i.e, $g_*$ is a strict tensor functor. Also, we require the inclusions $\iota_a \colon a \to m\otimes m$ to be respected by each $g_*$ in the sense that $g_*(\iota_a) = \iota_{g_*(a)},$ which can be safely assumed in our skeletal models of $\mathcal{TY}(A,\chi,\tau)$ and its automorphism group. 
\end{convention}

\section{Tambara-Yamagami Fixed Points} \label{sec:Theory}

\begin{proposition}\label{prop:monfun cocycle}
 With the notation above, if $\mu_g =1$ for all $g\in G$, then for all $g,h \in G$, and $a\in A$,  $J_{g,h}^a =1$ and $J^m_{-,-}$ is a 2-cocycle with coefficients in $\mathbb{Z}_2$.
\end{proposition} 

\begin{proof}
Consider the following commutative diagram:

\[\begin{tikzcd}[sep=small]
	{g_*(h_*(x\otimes y))} && {(gh)_*(x \otimes y)} \\
	\\
	{g_*(h_*(x)\otimes h_*(y))} && {(gh)_*(x) \otimes (gh)_*(y)} \\
	\\
	& {g_*(h_*(x)) \otimes g_*(h_*(y))}
	\arrow["{J_{g,h}^{x\otimes y}}", from=1-1, to=1-3]
	\arrow["{g_*(\mu_h^{x,y})}", from=3-1, to=1-1]
	\arrow["{\mu_{gh}^{x,y}}"', from=3-3, to=1-3]
	\arrow["{\mu_g^{h_*(x),h_*(y)}}", from=5-2, to=3-1]
	\arrow["{J^x_{g,h}\otimes J^y_{g,h}}"', from=5-2, to=3-3]
\end{tikzcd}\]
As $g_*$ is linear and $\mu_g=1$ for all $g\in G$, we get:

\[J_{g,h}^{x\otimes y}=J_{g,h}^x J_{g,h}^y\]
Applying the diagram when $x=m$ and $y=a \in A$, we see that:

$$J_{g,h}^mJ^a_{g,h} = J^{m \otimes a}_{g,h} = J^m_{g,h}$$
so $J^a_{g,h}=1$ for all $a\in A$. Furthermore, 

\[\begin{tikzcd}
	{g_1{_*}(g_2{_*}(g_3{_*}(m)))} && {(g_1g_2)_*({g_3}_*(m))} \\
	{g_1{_*}((g_2g_3)_*(m))} && {(g_1g_2g_3)_*(m)}
	\arrow["{J^{{g_3}_*(m)}_{g_1,g_2}}", from=1-1, to=1-3]
	\arrow["{{g_1}_*(J_{g_2,g_3}^m)}"', from=1-1, to=2-1]
	\arrow["{J_{g_1g_2,g_3}^m}", from=1-3, to=2-3]
	\arrow["{J_{g_1,g_2g_3}^m}"', from=2-1, to=2-3]
\end{tikzcd}\]
As $\mathcal{TY}$ has $\mathbb{Z}_2$ universal grading group with $m$ being the only nontrivially simple object, $g_*(m)=m$ for all $g\in G$ since functors preserve the universal grading. Therefore, the diagram implies:
$$J^m_{g_2,g_3}J^m_{g_1,g_2g_3} = J^m_{g_1,g_2}J^m_{g_1g_2,g_3}$$
so $J^m_{-,-}$ is a 2-cocycle. For each $a\in A$ we have the following diagram:

\[\begin{tikzcd}[column sep=3.5em, row sep=5em]
  g_*(h_*(a))
    \arrow[r, "{g_*(h_*(\iota_a))}"] 
    \arrow[d, "{J^a_{g,h}}"']
  & g_*(h_*(m\otimes m))
    \arrow[r, "{\scriptscriptstyle
      (\mu^{-1})_g^{h_*(m),h_*(m)}\circ g_*((\mu^{-1})_h^{m,m})
    }"{yshift=4pt}]
    \arrow[d, "{J^{m\otimes m}_{g,h}}"]
  & g_*(h_*(m))\otimes g_*(h_*(m))
    \arrow[d, "{J^m_{g,h}\otimes J^m_{g,h}}"]
  \\
  (gh)_*(a)
    \arrow[r, "{(gh)_*(\iota_a)}"']
  & (gh)_*(m\otimes m)
  \arrow[r, "{(\mu^{-1})_{gh}^{m,m}}"']
  & (gh)_*(m)\otimes (gh)_*(m)
\end{tikzcd}\]

Both squares individually commute, so the diagram commutes.  
Using the assumptions that $\mu_g=1$ (so that $\mu_g^{-1}=1$), the functor $g_*$ respects inclusions, and that $J^a_{g,h}=1$, we obtain $J^m_{g,h}\in \{\pm1\}$ as desired. 
\end{proof}

\begin{definition}[Adapted from \cite{Galindo2025}]\label{Def:G-inv_Sub_Cat} 
A full fusion subcategory $\mathcal{E} \subset \mathcal{TY}(A, \chi, \tau)$, is said to be $G$-\emph{invariant} if each $g_*$ restricts to an autoequivalence of $\mathcal{E}$.
\end{definition}

Next, we will consider a trivialization of a $G$-action on a proper $G$-invariant subcategory $\Vec{(B)} \subset \mathcal{TY}(A, \chi, \tau)$, where $B$ is a subgroup of $A$. The invariance property of $\Vec{(B)}$ provides an action of $G$ on $B$ by group homomorphisms.

The following definition describes how a trivialization of a $G$-action (Definition \ref{defn:Trivialization}) can be compatible with the conjugation action of $G$ on itself.  

\begin{definition}\label{Def:G-Triv} \cite{Galindo2025}
    Let $\mathcal{E}$ be a $G$-invariant subcategory of $\mathcal{C}$, and let $H\trianglelefteq G$. A trivialization $\eta$ of $F\mid_H$ on the subcategory $\mathcal{E}$ is $G$-\textit{equivariant} if the diagram below commutes:

\[\begin{tikzcd}
	{h_*(g_*(x))} && {g_*((g^{-1}hg)_*(x))} \\
	\\
	&& {g_*(x)}
	\arrow["can", from=1-1, to=1-3]
	\arrow["{\eta_h^{g_*(x)}}"', from=1-1, to=3-3]
	\arrow["{g_*(\eta^x_{g^{-1}hg})}", from=1-3, to=3-3]
\end{tikzcd}\]
Here $can$ is the morphism defined in section 3 of \cite{Galindo2025}.
\end{definition}

Written out explicitly, this diagram gives us the equation:

\begin{equation}\label{eqn:TrivializationEquation} g_*(\eta^x_{g^{-1}hg}J^x_{g^{-1}, hg}J^x_{h,g})(J^{h_*(g_*(x))}_{g,g^{-1}})^{-1} = \eta_h^{g_*(x)} \end{equation}
\noindent
where in the case $\mathcal{C} = \mathcal{TY}(A, \chi, \tau)$ and $\mathcal{E}=\Vec(B)$, by Proposition \ref{prop:monfun cocycle}, the equation simplifies to
\begin{equation}\label{eqn:EtaFixedPoint} \eta^x_{g^{-1}hg}=\eta_h^{g_*(x)}\end{equation}

\vspace{.5em}

Define a $G$-action $G \curvearrowright \Hom(H, \Hom(B, \mathbb{K}^\times))$ by $g \triangleright f(h)(b)= f(g \triangleright h)(g\triangleright b)= f(ghg^{-1})(g_*(b))$ where the left action on $B$ is inherited from how $G$ acts on $\Vec(B)$ and the left action on $H$ is by conjugation.

\begin{proposition}\label{prop: bijection}
Let $H\triangleleft G$ and $B$ be $G$-invariant and pointwise fixed by $H$. Let $G$ act on $\Hom(H, \Hom(B, \mathbb{K}^\times))$ as above. Then $G$-equivariant trivializations of $H$ on $\Vec(B)$ are in bijective correspondence with fixed points of $G \curvearrowright \Hom(H, \Hom(B, \mathbb{K}^\times))$.
\end{proposition} 

\begin{proof}
Let $\eta$ be a $G$-equivariant trivialization of $H$ on $\Vec(B)$. By definition, the diagrams 

    \[\begin{tikzcd}[column sep={1.5em, .1em, 1.5em}]
	{h_{1_*}(h_{2_*}(x))} & {h_{1_*}(x)} & {h_*(x)\otimes h_*(y)} & {h_*(x \otimes y)} \\
	{(h_1h_2)_*(x)} & x & {x\otimes y} & {x \otimes y}
	\arrow["{h_{1_*}(\eta^x_{h_2})}", from=1-1, to=1-2]
	\arrow["{J^x_{h_1,h_2}}"', from=1-1, to=2-1]
	\arrow["{\eta^x_{h_1}}", from=1-2, to=2-2]
	\arrow["{\mu_h^{x,y}}", from=1-3, to=1-4]
	\arrow["{\eta_h^x \otimes \eta_h^y}"', from=1-3, to=2-3]
	\arrow["{\eta^{x\otimes y}_h}", from=1-4, to=2-4]
	\arrow["{\eta^x_{h_1h_2}}", from=2-1, to=2-2]
	\arrow[equals, nfold, from=2-3, to=2-4]
\end{tikzcd}\]
are commutative. Since $\mu_h^{x,y}=J^x_{h_1,h_2}=1$ by assumption and Proposition \ref{prop:monfun cocycle} respectively, these diagrams translated into equations become:
\[\eta^x_{h_1h_2}  = \eta^x_{h_1}\eta^x_{h_2} \hspace{1em} \text{and} \hspace{1em} \eta^{x \otimes y}_h = \eta^x_h\eta^y_h\]
Define $f_\eta: H \times B \rightarrow \mathbb{K}^\times$ by $f_\eta(h,b) = \eta^b_h$. Then
\begin{align*}
f_\eta(h_1h_2,b) &= \eta^b_{h_1h_2}=\eta^b_{h_1}\eta^b_{h_2}=f_\eta(h_1,b)f_\eta(h_2,b) \\ f_\eta(h,b_1b_2) &= \eta^{b_1 b_2}_h = \eta_h^{b_1} \eta_h^{b_2} = f_\eta(h,b_1)f_\eta(h,b_2)
\end{align*}

By the second identity $f_\eta(h,-)$ is a character in $B$ for each $h$. This implies that $f_\eta\in \Hom(H, \Hom(B, \mathbb{K}^\times))$. By Equation \eqref{eqn:EtaFixedPoint}:
\[(g^{-1} \triangleright f_\eta)(h)(b) =f_\eta(g^{-1}hg)(g_*^{-1}(b)) =\eta^{g^{-1}_*(b)}_{g^{-1}hg} = \eta_h^{g_*(g^{-1}_*(b))}=\eta_h^b =f_\eta(h)(b)\]
Thus $f_\eta$ is fixed under the $G$-action.

Let $f \in \Hom(H, \Hom(B, \mathbb{K}^\times))$ be a fixed point of the $G$-action. As $B$ is pointwise fixed by $H$, so for $b\in B$, $h_*(b)= b$. Given each $b\in B$ and $h\in H$
 we produce a scalar $f(h)(b).$ Then we define $(\eta_f)_h^b:h_*(b)\rightarrow b$ by these scalars on simple objects. This defines a natural isomorphism between the functors $h_*$ and $\Id$.
 Moreover, $$\left(\eta_f\right)_h^{b_1}\left(\eta_f\right)_h^{b_2}=f(h)(b_1)f(h)(b_2) = f(h)(b_1b_2) =\left(\eta_f\right)_h^{b_1b_2} =\left(\eta_f\right)_h^{b_1 \otimes b_2}$$ 
since $f(h)\in \Hom(B, \mathbb{K}^\times)$ and $b_1b_2 \cong b_1\otimes b_2$. Moreover, since $\mu^{b_1,b_2}_h = \text{Triv}=1$, the following diagram commutes:
\[\begin{tikzcd}
	{h_*(b_1)\otimes h_*(b_2)} &&& {h_*(b_1\otimes b_2)} \\
	\\
	{\text{Id}(b_1)\otimes \text{Id}(b_2)} &&& {\text{Id}(b_1\otimes b_2)}
	\arrow["{\mu_h^{b_1,b_2}}", from=1-1, to=1-4]
	\arrow["{\left(\eta_f\right)^{b_1}_h \otimes \left(\eta_f\right)_h^{b_2}}", from=1-1, to=3-1]
	\arrow["{\left(\eta_f\right)_h^{b_1\otimes b_2}}", from=1-4, to=3-4]
	\arrow["{\text{Triv}}", from=3-1, to=3-4]
\end{tikzcd}\]

Thus $\left(\eta_f\right)_h: (h_*, \mu_h) \Rightarrow (\Triv, \Triv)$ is a monoidal natural transformation. Since the scalars are invertible $\left(\eta_f\right)_h$ is a monoidal natural isomorphism. It's clear to see that each step is invertible and thus these sets are in bijection. 

\end{proof}
Now consider the case where $\mathcal{E} = \mathcal{TY}(A, \chi, \tau)$. 

\begin{proposition}\label{prop: eta=1}
Let $\rho: G \rightarrow O(A,\chi)$ and $[J^m_{-,-}]\in H^2(G,\mathbb{Z}_2)$ parametrize the $G$-action on $\mathcal{TY}(A,\chi, \tau)$. Let $H\trianglelefteq G$, with $H\leq\ker\rho$. Then for such $H$, these trivializations are in bijection with $1$ cochains \(t:H\to\mathbb{Z}_2\), that satisfy $\delta t=J^m|_{H\times H}$ and the equation,
\[ t(g^{-1}hg)J^m_{h,g}=J^m_{g,g^{-1}hg}t(h),\hspace{3mm} g\in G, h\in H.\]
The trivialization corresponding to $t$ is determined by $\eta_a=\text{id}_a$ and $\eta_m=t(h)\text{id}_m$.
\end{proposition}

\begin{proof}
For $h \in H$, $\eta_h: (h_*, \mu_h) \Rightarrow (\Triv,\Triv)$ is a monoidal natural transformation, and so the following diagram commutes: 
\[\begin{tikzcd}
	{h_*(x)\otimes h_*(y)} && {h_*(x\otimes y)} \\
	\\
	{x\otimes y} && {x\otimes y}
	\arrow["{\mu_h^{x,y}}", from=1-1, to=1-3]
	\arrow["{\eta_h^x \otimes \eta_h^y}"', from=1-1, to=3-1]
	\arrow["{\eta^{x\otimes y}_h}", from=1-3, to=3-3]
	\arrow["{\text{Triv}}", from=3-1, to=3-3]
\end{tikzcd}\]
As $h_*$ is strict, $\eta^x_h\eta^y_h=\eta^{x\otimes y}_h$. Let $x=a$ and $y=m$, then as $a\otimes m =m$, $\eta^a_h\eta^m_h=\eta^m_h$, so $\eta^a_h=1$. Moreover, as $\eta: (F,J) \Rightarrow (\Triv, \Triv)$ is a monoidal natural transformation, the following diagram also commutes.

\[\begin{tikzcd}
	{F(h_1)\otimes F(h_2)} && {F(h_1\otimes h_2)} \\
	\\
	{\text{Tr}(h_1)\otimes \text{Tr}(h_2)} && {Tr(h_1\otimes h_2)}
	\arrow["{J_{h_1,h_2}}", from=1-1, to=1-3]
	\arrow["{\eta_{h_1} \otimes \eta_{h_2}}"', from=1-1, to=3-1]
	\arrow["{\eta_{h_1h_2}}"', from=1-3, to=3-3]
	\arrow["{\text{Triv}}", from=3-1, to=3-3]
\end{tikzcd}\]
Thus $J_{h_1,h_2}\eta_{h_1h_2} = \eta_{h_1}\eta_{h_2}$, so evaluating at $m$, we get $J_{h_1,h_2}^m\eta_{h_1h_2}^m = \eta_{h_1}^m\eta_{h_2}^m$. Therefore, $\delta(\eta^m)(h_1,h_2) = \eta_{h_1}^m\eta_{h_2}^m(\eta_{h_1h_2}^m)^{-1}=J^m_{h_1,h_2}$, so $\delta(\eta^m_{-}) = J^m_{-,-}$. Moreover, as $\eta^m_-$ is a $G$-invariant trivialization it satisfies \eqref{eqn:TrivializationEquation}. Let $k = g^{-1}hg$. Since $J^m_{-,-}$ is a normalized $2$-cocycle, so
\begin{small}\begin{align*}
J^m_{g,g^{-1}}&=J^m_{g,g^{-1}}J^m_{gg^{-1},g} = J^m_{g^{-1},g}J^m_{g,g^{-1}g}=J^m_{g^{-1},g}=J^m_{g^{-1},g}J^m_{g^{-1}g,k} \\ &=J^m_{g,k}J^m_{g^{-1},gk} = J^m_{g,g^{-1}hg}J^m_{g^{-1},hg}
\end{align*} \end{small}
so as $h_*(g_*(m))=m$
\[
g_*(\eta^m_{g^{-1}hg}J^m_{g^{-1}, hg}J^m_{h,g})(J^{h_*(g_*(m))}_{g,g^{-1}})^{-1} = \eta_h^{g_*(m)}
\]
becomes
\begin{equation} \label{eqn:CohomologyTrivializationRequirement}
\eta^m_{g^{-1}hg}J^m_{h,g} = J^m_{g,g^{-1}hg}\eta^m_h
\end{equation}

Conversely suppose that $t:H\rightarrow \mathbb{Z}_2$ is a 1-cochain satisfying $\delta t=J^m|_{H\times H}$ and
$t(g^{-1}hg)J^m_{h,g}=J^m_{g,g^{-1}hg}t(h),\hspace{3mm} g\in G, h\in H$. Define $\eta_h^a=\text{id}_a$ and $\eta_h^m=t(h)\text{id}_m$ and extend to direct sums linearly. This implies that $\eta_h$ is a natural isomorphism.  We note that these scalars assemble into a monoidal natural isomorphism. Since $\delta t=J|_{H\times H}$ it follows that $\eta$ will be compatible with multiplication in $H$, trivializing the $H$ action. Furthermore, as $t(g^{-1}hg)J^m_{h,g}=J^m_{g,g^{-1}hg}t(h)$ it follows that $\eta$ is a $G$-equivariant trivialization. This proves the bijection. 

\end{proof}

\begin{theorem} \label{thm:MainTheorem}
    Let $\rho \colon G \to O(A, \chi)$ be a homomorphism of groups, and $[\varphi] \in H^2(G, \mathbb{Z}/2\mathbb{Z})$ provide a $G$-action on $\mathcal{TY}(A, \chi, \tau)$. Then fusion subcategories of $\mathcal{TY}(A, \chi, \tau)^G$ occur in two families, which are respectively parameterized by 
    \begin{itemize}
        \item Triples $(H, B, \eta)$ such that 
        \begin{itemize}
            \item $H \trianglelefteq G$
            \item $B \leq A$ where $B$ is $G$-invariant and pointwise fixed by $H$
            \item $\eta \in \Hom(H, \Hom(B, \mathbb{K}^\times))$ a fixed point of the natural $G$ action. 
        \end{itemize}
        \item A normal subgroup $H \lhd G$ such that $H\leq \ker\rho$ with a trivialization $\eta$ of $[\varphi|_H]$ satisfying $$\eta(g^{-1}hg)\varphi(h,g) = \eta(h)\varphi(g,g^{-1}hg)$$
    \end{itemize}
\end{theorem}
\begin{proof} 
By Theorem \ref{thm:GJ subcats} fusion subcategories of $\mathcal{C}^G$ are parametrized by triples $(\mathcal{E},H,\eta)$. When $\mathcal{C}\cong \mathcal{TY}(A,\chi,\tau)$, there are two subcases: $\mathcal{E}$ a proper fusion subcategory and  $\mathcal{E}\cong \mathcal{TY}(A,\chi,\tau)$. In the first case, $\mathcal{E}\cong \Vec(B)$ and by Proposition \ref{prop: bijection} such triples are in bijection with the triples $(B,H,\eta)$ in the theorem statement. When $\mathcal{E}\cong \mathcal{TY}(A,\chi,\tau)$ the combination of Proposition \ref{prop:monfun G into Aut classification} and Proposition \ref{prop: eta=1} proves the second half of the theorem statement. 
\end{proof}

\begin{convention}
    Fusion subcategories of $\mathcal{TY}(A,\chi,\tau)^G$ corresponding to $(\Vec{(B)}, H, \eta)$ are said to be in the \textit{first family}, and those corresponding to $(\mathcal{TY}(A,\chi,\tau), H, \eta)$  in the \textit{second family}. 
\end{convention}

\section{Examples}\label{sec:Examples}

We make the results in the above section explicit for particular choices of $\mathbb{K},G,\rho, A, \chi$ and $\varphi$. 

Let $\mathbb{K} =\mathbb{C}$, $G = S_4$, $A = \mathbb{Z}_6 \times \mathbb{Z}_6$ and define the bicharacter $\chi$ such that 
$$\chi((x_1,y_1),(x_2,y_2)) = \exp\left( \frac{2\pi i}{6}(x_1y_2+x_2y_1) \right).$$ 
Consider $\phi, \psi \in O(\mathbb{Z}_6\times \mathbb{Z}_6, \chi)$ given by
\[
\phi(x,y) = (2x+3y, 3x+2y) \quad \& \quad \psi(x,y) = (x+3y, 3x+4y)
\]
Then $\phi^2 = \psi^3=\text{Id}$ and $\phi\circ\psi\circ \phi = \psi^{-1}$, so $\langle \phi, \psi\rangle \cong S_3$. Let $\gamma$ be that isomorphism from $S_3$ to $\langle \phi, \psi \rangle$ and let $\pi:S_4 \rightarrow S_3$ be the quotient map $S_4 \to S_3$. Define $\rho: S_4 \rightarrow O(\mathbb{Z}_6 \times \mathbb{Z}_6, \chi)$ by $\rho = \gamma\circ \pi$.
Let subgroups 
\[
\langle (2,0) \rangle, \langle (0,2) \rangle, \langle (2,2) \rangle, \langle (2,4) \rangle, \langle (2,0), (0,2) \rangle   
\]
be labeled $B_1$ through $B_5$ and the subgroup $\langle (3,0), (0,3)\rangle$ labeled as $K$. Let $KB_i$ denote the internal direct product. 

$H^2(S_4, \mathbb{Z}_2)$ has four classes. Since taking different representatives of the same class give an equivalent action, the number of trivializations for two representatives of the same class is the same.  Let 1 be the trivial cocycle and $\alpha, \beta, \alpha\beta$ the representatives of the other classes. Define
\[\alpha(g,h) = (-1)^{\epsilon(g)\epsilon(h)}\]
where $\epsilon:S_4 \rightarrow \{0,1\}$ is the sign homomorphism. The other two classes are difficult to distinguish - GAP was used to select a cocycle representative for a class which was neither trivial nor $\alpha$ and this class was declared to be $\beta$. This choice determines the other nontrivial class as $\alpha\beta$. 

\begin{theorem} \label{thm:Example}
Let $\chi$ be the hyperbolic bicharacter on $\mathbb{Z}_6\times \mathbb{Z}_6$. Let $(F,J): S_4  \rightarrow \text{Aut}_{\otimes}(\mathcal{TY}(\mathbb{Z}_6\times \mathbb{Z}_6, \chi, \tau))$ be classified by $\rho: S_4 \rightarrow O(\mathbb{Z}_6\times \mathbb{Z}_6, \chi)$ as given above and $[\varphi]\in H^2(S_4,\mathbb{Z}_2)$. Then $\mathcal{TY}(\mathbb{Z}_6\times \mathbb{Z}_6, \chi, \tau)^{S_4}$ has 53 fusion subcategories in the first family where the number of trivializations $\eta$ for each choice of $H, B$ are given by the tables:

\vspace{.25 em}

\begin{small} \begin{center} \label{Table1}
\begin{tabular}{cc}
\begin{tabular}{|c|c|c|c|c|}
\hline
$B$ & $e$ & $V_4$ & $A_4$ & $S_4$ \\
\hline
$e$ & 1 & 1 & 1& 1\\
\hline
$B_1$ &1&1&3&0 \\
\hline
$B_2$ &1&1&3&0 \\
\hline
$B_3$ &1&1&3&0 \\
\hline
$B_4$ &1&1&3&0 \\
\hline
$B_5$ &1&1&9&0 \\
\hline
\end{tabular} &
\begin{tabular}{|c|c|c|c|c|}
\hline
$B$ & $e$ & $V_4$ & $A_4$ & $S_4$ \\
\hline
$K$ & 1 & 2 & 0& 0\\
\hline
$KB_1$ &1&2&0&0 \\
\hline
$ KB_2$ &1&2&0&0 \\
\hline
$KB_3$ &1&2&0&0 \\
\hline
$KB_4$ &1&2&0&0 \\
\hline
$KB_5$ &1&2&0&0 \\
\hline
\end{tabular}
\end{tabular}
\end{center} \end{small}

\vspace{.25 em}

Moreover, the number of subcategories of the second family is 1 or 2 depending on the choice of $[\varphi]\in H^2(S_4, \mathbb{Z}_2)$. The number of trivializations $\eta$ per choice of $[\varphi]$ and $H$ are: \begin{center} \label{Table2}
\begin{tabular}{|c|c|c|c|c|}
\hline
& 1 & $\alpha$ & $\beta$ & $\alpha\beta$ \\
\hline 
$e$ & 1& 1& 1 &1 \\
\hline
$V_4$ & 1 & 1& 0& 0\\
\hline
\end{tabular}
\end{center}
\end{theorem}
\subsection{Proof of Theorem \texorpdfstring{\ref{thm:Example}}{4.1}}

\begin{definition}[\cite{MR1324339}]
Let $G$ be a group and $M$ a $G$-module, then the group of co-invariants of $M$, denoted $M_G$, is defined to be the quotient $M / \langle g\cdot m-m\rangle $.\end{definition}

\begin{lemma} \label{lemma:EffectiveCohomology}
Let $G$ be a group and $B$ a finite abelian group. Let $N$ be a finite group such that $G$ acts on $B$ and $N$ by automorphisms. Then $$|\Hom(N, \widehat{B})^G| = |(N_{ab}\otimes B)_G|$$ 

\end{lemma}

\begin{proof}
As $\widehat{B}$ is abelian, by the tensor-hom adjunction:
\[
\Hom(N, \widehat{B}) \cong  \Hom(N_{ab}, \Hom(B, \mathbb{C}^\times)) \cong \Hom(N_{ab} \otimes B, \mathbb{C}^\times ).
\]
Here $G$ acts on $D=N_{ab} \otimes B$ via the diagonal action.
Let $f \in \widehat{D}^G$, then for all $g\in G$ and $a\in N_{ab}\otimes B$, $f(g\cdot a)=f(a)$, so $f(g \cdot a-a)=1$. Therefore, $f$ is trivial on the subgroup generated by $\{g\cdot a-a| a\in N_{ab}\otimes B, g\in G\}$. Therefore, $f$ factors through $(N_{ab}\otimes B)/\langle g\cdot a-a \rangle = (N_{ab} \otimes B)_G$, and since $f$ was arbitrary, $\widehat{D}^G \cong \Hom((N_{ab}\otimes B)_G,\mathbb{C}^\times )$. Finally, as $N$ is finite,  $(N_{ab}\otimes B)_G$ is finite abelian, so $(N_{ab}\otimes B)_G \cong \Hom((N_{ab}\otimes B)_G,\mathbb{C}^\times )$. Thus
\[
\Hom(N, \Hom(B,\mathbb{C}^\times))^G \cong \widehat{D}^G \cong D_G
\] and so $|\Hom(N, \Hom(B,\mathbb{C}^\times))^G| = |(N_{ab}\otimes B)_G|$. 

\end{proof}

As $\rho: S_4 \rightarrow O(\mathbb{Z}_6\times \mathbb{Z}_6, \chi)$ is defined by $\rho =\gamma \circ \pi$, we obtain $\rho((123)) = \psi$ and $\rho((12)) = \phi$. \\
Next, we examine the $G$-invariant subgroups of $\mathbb{Z}_6\times \mathbb{Z}_6$. Notice that 
$$\mathbb{Z}_6 \times \mathbb{Z}_6 \cong (\mathbb{Z}_3 \times \mathbb{Z}_2)\times (\mathbb{Z}_3\times \mathbb{Z}_2)\cong (\mathbb{Z}_3\times \mathbb{Z}_3)\times (\mathbb{Z}_2
\times \mathbb{Z}_2).$$
Therefore, every subgroup of $\mathbb{Z}_6 \times \mathbb{Z}_6$ is a direct product of a subgroup of the $\mathbb{Z}_2\times \mathbb{Z}_2$ component and the $\mathbb{Z}_3\times \mathbb{Z}_3$ component. As $\mathbb{Z}_2\times \mathbb{Z}_2$ has 5 subgroups and $\mathbb{Z}_3 \times \mathbb{Z}_3$ has 6 subgroups, $\mathbb{Z}_6 \times \mathbb{Z}_6$ has 30 subgroups. Of these, not all will be $G$-invariant. Since
\begin{align*}
(123)\triangleright (3,0) &= \rho(123)((3,0)) = \psi((3,0)) = (3,3)\\
(123)\triangleright (3,3) &= \rho(123)((3,3)) = \psi((3,3)) = (0,3)\\
(123)\triangleright (0,3) &= \rho(123)((0,3)) = \psi((0,3)) = (3,0)
\end{align*}
any $G$-invariant subgroup containing any of the elements $(0,3), (3,0)$, or $(3,3)$ contains all of them. Therefore, the only subgroups of the $\mathbb{Z}_2 \times \mathbb{Z}_2$ component that are $G$-invariant are the trivial subgroup and the entire subgroup. Furthermore, as 
\begin{align*}
    \phi(2,0) &=(4,0) \quad \quad \quad  \psi(2,0) = (2,0) \\
    \phi(0,2) & = (0,4) \quad \quad \quad \psi(0,2) = (0,2)
\end{align*}
so $\phi(2,0),\psi(2,0) \in \langle(2,0) \rangle $ and $\phi(0,2),\psi(0,2) \in \langle (0,2) \rangle$. Thus the G-action fixes the subgroups of the $\mathbb{Z}_3 \times \mathbb{Z}_3$ component, and so they are all $G$-invariant. Thus in total $\mathbb{Z}_6 \times \mathbb{Z}_6$ has 12 $G$-invariant subgroups, given by $e, B_i, K$ and $KB_i$ for $1\le i\le 5$. With this we can calculate the number of fusion subcategories of $\mathcal{TY}(A,\chi,\tau)^{S_4}$ belonging to the first family for each choice of $H$. 

\subsubsection{The trivial subgroup}  Since all subgroups are pointwise fixed by the trivial subgroup, for any subgroup $B$, $\Hom(e, \Hom(B, \mathbb{C}^\times))$ consists only of the trivial homomorphism, and since the trivial homomorphism is always fixed, $|\Hom(e, \Hom(B, \mathbb{C}^\times))^{S_4}|=1$.

\subsubsection{\texorpdfstring{The subgroup $V_4$}{Klein four subgroup}}
Now for $H = V_4$, since $V_4$ is the kernel of $\rho$, all subgroups are pointwise fixed by $V_4$. As the orders of $V_4$ and $B_i$ are co-prime, $|\Hom(V_4,\Hom(B_i, \mathbb{C}^\times))|=1$, so 
\[\left\lvert\Hom(V_4,\Hom(B_i, \mathbb{C}^\times))^{S_4}\right\lvert=1.\]
$V_4 \otimes K$ has order sixteen and is generated by $$\{(12)(34)\otimes (3,0), (12)(34)\otimes (0,3), (13)(24)\otimes (3,0), (13)(24)\otimes (0,3) \}.$$ We can compute:
\begin{scriptsize}
\begin{align*}
&(12) \cdot((12)(34)\otimes (3,0))-(12)(34)\otimes (3,0) = (12)(34)\otimes (3,3) \\
& (13)\cdot ((13)(24)\otimes (0,3))-(13)(24)\otimes (0,3) = (13)(24)\otimes (3,0)\\
& (23)\cdot ((14)(23)\otimes (3,0))-(14)(23)\otimes(3,0) = (14)(23)\otimes (0,3)\\
&(132)\cdot ((12)(34)\otimes (3,0))-(12)(34)\otimes(3,0) = (13)(24)\otimes(0,3)+(12)(34)\otimes(3,0)
\end{align*}
\end{scriptsize}
By combining these equations, 8 distinct elements of $V_4 \otimes K$ can be generated, so $|\langle g\cdot a-a \rangle |\ge 8$. Moreover, calculating for each $g\in S_4$ and $a\in V_4\otimes K$, it can be shown $g\cdot a-a$ always takes one of the previously calculated eight values, so $|\langle g\cdot a-a\rangle|=8$. Therefore,  $|(V_4\otimes K)_{S_4}|=2$ and by Lemma \ref{lemma:EffectiveCohomology},
\[\left\lvert\Hom(V_4,\Hom(K, \mathbb{C}^\times))^{S_4}\right\rvert=2\]

Now consider $KB_i$. Since every element in $KB_i$ can be written as $kb$ with $k\in K$, $b\in B_i$ and $K\cap B_i = \{0\}$, a homomorphism $f\in \Hom(V_4, \Hom(K, \mathbb{C}^\times))^{S_4}$ can be extended to a homomorphism $\tilde{f}\in \Hom(V_4, \Hom(KB_i, \mathbb{C}^\times))$ by defining $\tilde{f}(h)(k) = f(h)(k)$ for $k\in K$ and $\tilde{f}(h)(b)=1$ for $b\in B_i$ for all $h\in V_4$. Moreover, as 
 \begin{small}\begin{align*}
 &g\cdot \tilde{f}(h)(kb)  =\tilde{f}(ghg^{-1})(\rho(g)(kb)) = \tilde{f}(ghg^{-1})(\rho(g)(k))\tilde{f}(ghg^{-1})(\rho(g)(b))\\ &= (f(ghg^{-1})(\rho(g)(k)))(1) = f(h)(k) = \tilde{f}(h)(k) = (\tilde{f}(h)(k))(\tilde{f}(h)(b))\\ &=\tilde{f}(h)(kb)
 \end{align*} \end{small}
 since $K$ and $B_i$ are $G$-invariant, so $\rho(g)(k)\in K$ and $\rho(g)(b)\in B_i$. Thus $\tilde{f}$ is fixed under the $G$-action and thus $|\Hom(V_4, \Hom(K, \mathbb{C}^\times))^{S_4}|\le |\Hom(V_4, \Hom(KB_i, \mathbb{C}^\times))^{S_4}|$. $V_4\otimes KB_i \cong (V_4 \otimes K)\oplus (V_4 \otimes B_i) \cong V_4\otimes K$ due to the orders of $V_4$ and $B_i$ being co-prime. Thus $|V_4 \otimes KB_i|=16$. For $a\in V_4 \otimes K$, then $a\in V_4 \otimes KB_i$, we get $|\langle g\cdot a-a \rangle |\ge 8$ by previous calculation. Thus $|(V_4 \otimes KB_i)_{S_4}| =16/|\langle g \cdot a-a\rangle |\le 2$ and so by Lemma \ref{lemma:EffectiveCohomology}
\[
\left\lvert \Hom(V_4, \Hom(KB_i, \mathbb{C}^\times))^{S_4}\right\rvert =2
\]

\subsubsection{\texorpdfstring{The subgroup $A_4$}{Alternating group}} Now we address $H = A_4$. Since $\rho(A_4) = \langle \psi \rangle$, the trivial subgroup and subgroups $B_i$ are pointwise fixed by $A_4$. However, $\psi$ permutes the elements of $K$, so the subgroups $K$ and $KB_i$ for $1\le i\le 5$ are not pointwise fixed by $A_4$. Therefore, there are no triples of the form $(A_4, K, \eta)$ nor of the form $(A_4, KB_i, \eta)$ for $1\le i\le 5$. To calculate the number of trivializations for subgroups $B$ pointwise fixed by $A_4$, first consider the order 3 subgroups, $B_i$ for $i\le 4$. 
By Lemma \ref{lemma:EffectiveCohomology},
\begin{small}\[\left\lvert\Hom(A_4, \Hom(B_i, \mathbb{C}^\times))^{S_4}\right\rvert = \left\lvert((A_4)_{ab}\otimes B_i)_{S_4}\right\rvert\] \end{small}
Since $(A_4)_{ab}\otimes B_i \cong \mathbb{Z}_3$, $|((A_4)_{ab}\otimes B_i)_{S_4}| = 3/|\langle g\cdot a-a \rangle |$. 
Let $e, (123), (132)$ be the representatives for the elements of $(A_4)_{ab}$. For $g$ odd:
\begin{align*} &g \cdot ((123) \otimes (0,2)) = g\cdot (123) \otimes g\cdot (0,2) = g(123)g^{-1}\otimes \phi(0,2) = \\&(132)\otimes (0,4)= (132)\otimes 2(0,2) = (132)^2 \otimes (0,2) = (123)\otimes (0,2) \end{align*}
so $g \cdot ((123)\otimes (0,2))-(123)\otimes (0,2)=0$. Repeating the previous calculation for the other elements, we get $g\cdot a-a=0$ for all $a \in (A_4)_{ab}\otimes B_i$ and all $g \in S_4$. Thus $\langle g\cdot a-a\rangle  = \{0\}$, so  $|((A_4)_{ab}\otimes B_i)_{S_4}|=3$, and thus 
\[\left\lvert\Hom(A_4, \Hom(B_i, \mathbb{C}^\times))^{S_4}\right\rvert=3\]

Now consider $B_5$. As $B_5 \cong \mathbb{Z}_3 \times \mathbb{Z}_3$, $(A_4)_{ab}\otimes B_5 \cong \mathbb{Z}_3\oplus \mathbb{Z}_3$, so $|(A_4)_{ab}\otimes B_5|=9$. Let $h \in (A_4)_{ab}$ and $b \in B_5$. Then $b\in B_i$ for some $i\le 4$. Therefore, by above calculation $g\cdot (h\otimes b)-h\otimes b=0$. Thus $\langle g \cdot (h\otimes b)-(h\otimes b) \rangle =\{0\}$, so $|((A_4)_{ab}\otimes B_5)_{S_4}|= 9$. Therefore, by Lemma \ref{lemma:EffectiveCohomology}:  \[\left\lvert\Hom(A_4, \Hom(B_5,\mathbb{C}^\times))^{S_4}\right\rvert = 9\]

\subsubsection{The entire group} For $H=S_4$, since 
\begin{align*}
    (12) & \triangleright (2,0) = \rho(12)(2,0) = \phi(2,0)=(4,0)\\
    (12) & \triangleright (0,2) = \rho(12)(0,2) = \phi(0,2)=(0,4)\\
    (123) & \triangleright (3,0) = \rho(123)(3,0)=\psi(3,0)=(3,3)
\end{align*}
Only the trivial subgroup is fixed pointwise by $S_4$. Thus for $B\neq e$, there are no triples $(S_4, B, \eta)$ and by the same argument used for the trivial subgroup
\[
\left\lvert\Hom(S_4, \Hom(e, \mathbb{C}^\times))^{S_4}\right\rvert =1
\]

Thus there are $53$ fusion subcategories of $\mathcal{TY}(\mathbb{Z}_6\times \mathbb{Z}_6, \chi, \tau)^{S_4}$ belonging to the first family. 

\subsubsection{\texorpdfstring{Cohomology}{Cohomology}}
The number of fusion subcategories in the second family depends on the cohomology class of $J^m_{-,-}$. We calculate the results for each of the four classes in $H^2(S_4, \mathbb{Z}_2)$. Since $\ker(\rho)=V_4$, trivializations only exist when $H=e$ or $H=V_4$. When $H=e$, since the representatives are normalized, $\eta^m_e=1$ which satisfies Equation \eqref{eqn:TrivializationEquation}. Thus for $H=e$, each representative has exactly one trivialization. 
For $H=V_4$, $\alpha$ restricted to $V_4$ is the same as the trivial representative, taking the value 1 on all pairs. Thus a trivialization, $\eta^m_-$, for either representative $1$ or $\alpha$  satisfies
\[
\eta^m_{h_1}\eta^m_{h_2}=\eta_{h_1h_2}^m
\]
so $\eta^m_-\in \Hom(V_4, \mathbb{Z}_2)$. Moreover, for the trivial representative and $\alpha$, equation \eqref{eqn:CohomologyTrivializationRequirement} simplifies to equation \eqref{eqn:EtaFixedPoint} since $\alpha(h,g)=\alpha(g,g^{-1}hg)=1$ for all $h\in V_4$ and $g\in S_4$. Therefore,
\[
\eta^m_{g^{-1}hg}=\eta^m_h
\]
for $g\in S_4$. As the trivial character on $V_4$ is the only one fixed under conjugation by $S_4$, for representatives $1$ and $\alpha$ there is exactly one trivialization. 
For both $\beta$ and $\alpha\beta$, using representatives obtained using GAP \cite{GAP4.12.1},  $\beta((12)(34), (12)(34))= \alpha\beta((12)(34),(12)(34))=-1$. Assume that a $V_4$ trivialization $\eta$ exists. Since $\beta(e,e)=\alpha\beta(e,e)=1$, we get that $\eta_e^m=1$. Thus 
\[
\eta^m_{(12)(34)}\eta^m_{(12)(34)}=\eta^m_eJ^m_{(12)(34),(12)(34)} =-1
\]
so $(\eta^m_{(12)(34)})^2=-1$ while $\eta^m_{(12)(34)}\in \{\pm 1\}$. Thus we have reached a contradiction and $V_4$ has no trivializations for representatives $\beta$ and $\alpha\beta$. 

\bibliographystyle{alpha}
\bibliography{bibliography}

\end{document}